\documentstyle[12pt]{article}

\catcode`\@=11

\def\section{\@startsection {section}{1}{\z@}{-3.5ex plus 
  -1ex minus -.2ex}{2.3ex plus .2ex}{\normalsize\bf}}
\def\subsection{\@startsection {subsection}{2}{\z@}{-3.25ex plus 
  -1ex minus -.2ex}{1.5ex plus .2ex}{\normalsize\bf}}

\catcode`\@=12

\newtheorem{thm}{Theorem}[section]
\newtheorem{prop}[thm]{Proposition}
\newtheorem{cor}[thm]{Corollary}
\newtheorem{lemma}[thm]{Lemma}

\newcommand{\Fp}{{{\bf F}_{\kern-.125em p}}}
\newcommand{\Fq}{{{\bf F}_{\kern-.125em q}}}
\newcommand{\Fpb}{{{\overline{\bf F}}_{\kern-.125em p}}}
\newcommand{\Fqb}{{{\overline{\bf F}}_{\kern-.125em q}}}
\newcommand{\Kb}{{\overline{K}}}
\newcommand{\nb}{{\gcd(n,q-1)}}
\newcommand{\Xb}{{\overline{X}}}
\newcommand{\sigmab}{{\overline{\sigma}}}
\newcommand{\Mb}{{\overline{\cal M}}}
\newcommand{\wh}{{\hat{w}}}
\newcommand{\vh}{{\hat{v}}}
\newcommand{\BA}{{\bf A}}
\newcommand{\BE}{{\bf E}}
\newcommand{\BP}{{\bf P}}

\newcommand{\BZ}{{\bf Z}}
\newcommand{\hollowsquare}{\vbox{\hrule\hbox{\vrule\kern.6em\vbox{\kern.6em}\vrule}\hrule}}
\newcommand{\ra}{\rightarrow}
\newcommand{\lra}{\longrightarrow}
\newcommand{\abs}{{\rm abs}}
\newcommand{\Aut}{\mathop{\rm Aut}\nolimits}
\newcommand{\can}{{\mathop{\rm can}\nolimits}}
\newcommand{\Gal}{\mathop{\rm Gal}\nolimits}
\newcommand{\Ig}{\mathop{\rm Ig}\nolimits}
\newcommand{\ord}{\mathop{\rm ord}\nolimits}
\newcommand{\Spec}{\mathop{\rm Spec}\nolimits}
\newcommand{\SL}{\mathop{\rm SL}\nolimits}
\newcommand{\mapright}[1]{\mathop{\lra}\limits^{#1}}

\newcommand{\version}{19920303}

\newenvironment{proof}{\begin{sc}\noindent Proof: \end{sc}}{
     \hbox to 2em{}\nobreak\hfill$\hollowsquare$\par\medskip}
\newenvironment{notation}{\noindent\begin{bf}Notation: \end{bf}}{}
\newenvironment{defn}{\noindent\begin{bf}Definition: \end{bf}}{}
\newenvironment{remark}{\noindent\begin{bf}Remark: \end{bf}}{}

\begin{document}

\noindent
\begin{center}
\scriptsize \bf On the group orders of elliptic curves over finite fields
\rm --- Everett W.  Howe
\hfill version \version
\end{center}
\hrule

\vskip0.8cm
\noindent
\bf
On the group orders of elliptic curves over finite fields

\noindent
\rm ({\it Compositio Mathematica\/} {\bf 85} (1993) 229--247)

\medskip\noindent\rm
Everett W. Howe

\bigskip\noindent
\baselineskip=10pt
\footnotesize
Department of Mathematics\hfil\break
University of California\hfil\break
Berkeley, CA 94720 

\bigskip\noindent
\bf Abstract.
\rm Given a prime power $q$, for every pair of positive integers $m$ and $n$ 
with $m | \gcd(n, q-1)$ we construct a modular curve over 
${\scriptstyle\bf F}_{\kern-.125em q}$
that parametrizes elliptic curves over 
${\scriptstyle\bf F}_{\kern-.125em q}$ along
with ${\scriptstyle\bf F}_{\kern-.125em q}$-defined points $P$ and $Q$ of 
order $m$ and~$n$, respectively, with
$P$ and ${n \over m}Q$ having a given Weil pairing.  Using these curves,
we estimate the number of elliptic curves over 
${\scriptstyle\bf F}_{\kern-.125em q}$ that have
a given integer $N$ dividing the number of their 
${\scriptstyle\bf F}_{\kern-.125em q}$-defined points.


%
% 11         -- Number theory
%   11G      -- Arithmetic algebraic geometry (Diophantine geometry)
%      11G20 -- Curves over finite and local fields
% 14         -- Algebraic geometry
%   14G      -- Arithmetic problems.  Diophantine geometry
%      14G15 -- Finite ground fields
%   14H      -- Curves
%      14H52 -- Elliptic curves

\smallskip\noindent
\bf Acknowledgements.
\rm The author was supported by a United States Department of Education
National Need Fellowship.  The author would like to thank 
Hendrik Lenstra for suggesting the proof of Proposition~\ref{form_sum:prop},
and Bas Edixhoven and Hendrik Lenstra for providing  helpful comments
and advice.

\normalsize
\baselineskip=18pt

\section{Introduction}

Given a prime power $q$ and a positive integer $N$, it is natural to wonder
how likely it is for a randomly chosen elliptic curve over $\Fq$ to have $N$ 
dividing the number of its $\Fq$-defined points.  The purpose of this paper is
to make sense of this question and to provide an estimate for its answer.

Since $\Fq$-isomorphic curves have the same number of $\Fq$-defined points,
we will only be interested in $\Fq$-isomorphism classes of elliptic 
curves over $\Fq$.  In particular, we will look at the set
\[ V(\Fq;N) = \Big\{ E/\Fq : N | \#E(\Fq) \Big\}\Big/\cong_\Fq ;\]
we want to know how large this set is, compared to the set of all
$\Fq$-isomorphism classes of elliptic curves over $\Fq$.  However,
it will be easiest to estimate not the usual cardinality of $V(\Fq;N)$ but 
rather the {\em weighted cardinality} of $V(\Fq;N)$, where the {\em weighted 
cardinality} of a set $S$ of $\Fq$-isomorphism classes of elliptic curves
over $\Fq$ is defined to be
\[\#'S = \sum_{[E]\in S} {1\over\#\Aut_\Fq(E)} ,\]
where $[E]$ denotes the $\Fq$-isomorphism class of the elliptic curve $E$.
Often, formulas for weighted cardinalities of such sets $S$ work out 
better than formulas for the usual cardinalities; for instance, we will
see in Corollary \ref{num_curves:cor} that 
\begin{equation}
\label{num_curves:eqn}
\#' \Big\{ E : \hbox{ $E$ is an elliptic curve over $\Fq$}\Big\} 
                                              \Big/\cong_\Fq \ = q ,
\end{equation}
whereas the corresponding formula for the ordinary cardinality
depends on the value of $q \bmod 12$.  In any case, since 
$\Aut_\Fq(E)=\{\pm 1\}$ except possibly when $j(E)$ is $0$ or $1728$
(see \cite{Si}, section III.10), the weighted cardinality of such a set 
$S$ is generally about half of its usual cardinality.

In view of (\ref{num_curves:eqn}), we will interpret
the ratio $\#'V(\Fq;N) / q $  as the probability
that a random elliptic curve over $\Fq$ has $N$ dividing the number
of its $\Fq$-defined points.  The following theorem gives an
estimate of this ratio.

\begin{thm}
\label{theorem:thm}
There is a constant $C\le1/12+5\sqrt{2}/6\approx 1.262$
such that the following statement is true:
Given a prime power $q$, let $r$ be the multiplicative arithmetic function 
such that for all primes $\ell$ and positive integers $a$
\[r(\ell^a)=\left\{ \begin{array}{ll}
            \displaystyle {1\over \ell^{a-1}(\ell-1)}&
                         \mbox{if $q\not\equiv 1 \bmod{\ell^c};$} \\
     \ & \\
            \displaystyle{{\ell^{b+1}+\ell^b-1}\over{\ell^{a+b-1}(\ell^2-1)}}&
                         \mbox{if $q\equiv 1 \bmod{\ell^c},$} 
                            \end{array}
                     \right.   \]
where $b=\lfloor{a/2}\rfloor$, the greatest integer less than or equal 
to $a/2$,
and $c=\lceil{a/2}\rceil$, the least integer greater than or equal
to $a/2$.
Then for all positive integers $N$ we have
\begin{equation}
\label{prob_est:eqn}
\left|{\#'V(\Fq;N)\over q}-r(N)\right| 
                      \le{CN\rho(N)2^{\nu(N)}\over\sqrt{q}}\ ,
\end{equation}
where $\rho(N)=\prod_{p| N}((p+1)/(p-1))$ and 
$\nu(N)$ denotes the number of prime divisors of $N$.
\end{thm}

It is interesting to note that $r(N)$ is greater than $1/N$ and for
many values of $N$ is not much less than $1/\varphi(N)$.  Thus, loosely
speaking, when $q$ is large with respect to $N$ it is more likely that 
a random elliptic curve over $\Fq$ has 
$N$ dividing its number of points than it is that a random integer is
divisible by $N$.

H.W. Lenstra, Jr.~has proven the inequality~(\ref{prob_est:eqn}) in the 
special case when $N$ and $q$ are distinct primes with $q>3$ (see \cite{Le}, 
Proposition 1.14, page 660).  
Lenstra's proof depends on properties of modular curves over $\Fp$; in 
particular, he uses the modular curves $X(\ell)$ and $X_1(\ell)$, for 
primes $\ell\ne p$.  My extension of Lenstra's proposition is obtained
by extending his proof, and accordingly my proof will require the 
study of modular curves which I will denote $X_q(m,n)$.

In section~\ref{form:sec}, I briefly prove some results about forms that
will be needed in sections~\ref{X:sec} and~\ref{proof:sec}. In 
section~\ref{X:sec}, I define the curves $X_q(m,n)$ as quotients of 
more familiar modular curves, give a modular interpretation of their
$\Fq$-defined points, and use Weil's estimate to approximate the
number of their $\Fq$-defined points.  Finally, in section~\ref{proof:sec}
I use the interpretation and bounds of section~\ref{X:sec}
for a number of curves to prove Theorem \ref{theorem:thm}.

\begin{notation}
Throughout this paper, if $C$ is a curve over a field $K$, and
if $L$ is an extension field of $K$, we will denote by $C_L$ the $L$-scheme
$C \times_{\Spec(K)} \Spec(L)$.  Similarly, if $P$ is a $K$-defined point on 
such a curve $C$, we will denote by $P_L$ the point on $C_L$ obtained from 
$P$ by base extension.  
If $E$ is an elliptic curve over $K$ with zero point $O$, then the curve $E_L$ 
has a unique structure of an
elliptic curve over $L$ with zero point $O_L$; 
when we mention the curve $E_L$, we will be referring to it as an elliptic
curve, unless we explicitly state otherwise.
The letters $p$ and $\ell$ are reserved for prime numbers. 
For real numbers $x$, we will denote by $\lfloor x\rfloor$ the greatest
integer less than or equal to $x$ and by $\lceil x\rceil$ the least integer
greater than or equal to $x$.
Also, we will make use of five arithmetic functions:
the M\"obius function~$\mu$;  the function $\nu$ such that $\nu(n)$ 
is the number of prime divisors of $n$; the Euler totient function $\varphi$, 
defined by $\varphi(n)=n\prod_{p|n}(1-1/p)$; the function $\psi$ defined
by $\psi(n)=n\prod_{p|n}(1+1/p)$; and the function $\rho$ defined by
$\rho(n)=\prod_{p|n}((p+1)/(p-1))$.
\end{notation}

\section{Forms}
\label{form:sec}

\begin{defn}
Let $E$ be an elliptic curve over a field $K$, and let $L$ be an extension
field of $K$. An elliptic curve $E'$ over $K$ is called an
{\em $L/K$-form of $E$} (or simply a {\em form of $E$}, if $L$ and $K$ are
clear from context) if $E_L$ and $E'_L$ are isomorphic over $L$.
We denote by $\BE(L/K;E)$ or simply $\BE(E)$ the set of
forms of $E$, up to $K$-isomorphism:
\[\BE(L/K;E) = \Big\{E'/K : E'_L \cong E_L\Big\}\Big/{\cong_K} ; \]
and we denote by $[E']_K$  or simply $[E']$ the $K$-isomorphism class of $E'$.
Suppose we are also given points $P,Q\in E(K)$.  A triple $(E',P',Q')$,
where $E'$ is an elliptic curve over $K$ and $P'$ and $Q'$ are
points of $E'(K)$, is 
called an {\em $L/K$-form of $(E,P,Q)$} if there is an $L$-isomorphism from
$E_L$ to $E'_L$ that takes $P_L$ to $P'_L$ and $Q_L$ to $Q'_L$.  We denote 
by $\BE(E,P,Q)=\BE(L/K;E,P,Q)$ the set of $L/K$-forms of $(E,P,Q)$, up to 
$K$-isomorphism, and we denote by $[E',P',Q']_K$ the $K$-isomorphism class
of the triple $(E',P',Q')$.
\end{defn}  

Suppose $L$ is a finite or infinite Galois extension of $K$ with topological
Galois group $G$, and suppose $E$ is an elliptic curve over $K$.
Let $A$ be the finite group $\Aut_L(E_L)$ of all $L$-automorphisms of 
$E_L$, and let $B$ be the group of all commutative diagrams
% Let's draw a commutative diagram.
\begin{center}
\begin{picture}(700,700)(-350,-350)
%
% First, the four schemes...
\put(-350, 250){\makebox(0,0){$E_L$}}
\put( 350, 250){\makebox(0,0){$E_L$}}
\put(-350,-250){\makebox(0,0){$\Spec(L)$}}
\put( 350,-250){\makebox(0,0){$\Spec(L)$}}
%
% Next, the two vertical arrows...
\put(-350, 180){\vector( 0,-1){360}}
\put( 350, 180){\vector( 0,-1){360}}
%
% Then the two horizontal arrows...
\put(-270, 250){\vector( 1, 0){540}}
\put(-150,-250){\vector( 1, 0){300}}
%
% And finally, the names of the horizontal morphisms.
\put(   0, 300){\makebox(0,0){$\alpha$}}
\put(   0,-200){\makebox(0,0){$\sigmab$}}
\end{picture}
\end{center}
where $\alpha$ is an automorphism of $E_L$ as a $K$-scheme that fixes the
zero point of $E_L$,
and where for any element $\sigma$ of $G$ we denote by $\sigmab$ the 
scheme automorphism of $\Spec(L)$ obtained from the field automorphism
$\sigma^{-1}$ of $L$.
There is clearly an exact sequence of groups
\begin{equation}
\label{seq:eqn}
 1 \lra A \lra B \mapright{\pi} G \lra 1 
\end{equation}
where $\pi$ is the projection map taking an element $(\alpha, \sigmab)$
of $B$ to the element $\sigma$ of $G$.  The sequence (\ref{seq:eqn}) has
a canonical splitting $G\ra B$ defined by sending $\sigma\in G$ to
the element $(1\times\sigmab, \sigmab)$ of $B$, where 
$1\times\sigmab$ is the $K$-scheme automorphism of 
$E_L=E\times_{\Spec(K)} \Spec(L)$ 
obtained by fixing $E$ and applying $\sigmab$ to $\Spec(L)$.
As a set, $B$ is the product of $A$ and $G$; if we give $A$ the discrete
topology and $B$ the product topology, the sequence (\ref{seq:eqn}) is
even an exact sequence of topological groups.

{}From \cite{Se:CG} (see in particular section III.1.3), we know that 
$\BE(L/K;E)$ is isomorphic (as a set with
a distinguished element) to the cohomology set $H^1(G,A)$, where the
cohomology is in the sense of section I.5 of~\cite{Se:CG} (see also~\cite{Si},
sections X.2 and X.5).  A cocycle, in this sense,
corresponds to a continuous homomorphism $s:G\ra B$ splitting the 
exact sequence~(\ref{seq:eqn}); such a section gives an action of 
$G$ on $E_L$, and this defines by Galois descent an elliptic curve $E(s)/K$ 
and an isomorphism $f_s:E_L\ra E(s)_L$, unique up to $\Aut_K(E(s))$
--- see~\cite{We:FD} or section~V.20 of~\cite{Se:AG} for the case of
finite extensions $L/K$, and
compare problem~II.4.7 (page~106) of \cite{Ha}.  The group $A$ acts on
the set $S$ of sections by conjugation, and two cocycles are cohomologous if
and only if  their 
associated sections lie in the same $A$-orbit of $S$.  Also, the stabilizer
of a section $s$ is isomorphic to the group of $K$-automorphisms of the
associated form $E(s)$.  Thus the orbit-decomposition formula
(\cite{La}, page 23) gives 
\begin{equation}
\label{gen_form_sum:eqn}
\sum_{[E']\in\BE(E)}{{\#A}\over\#\Aut_K(E')} = \#S .
\end{equation}

\begin{prop}
\label{form_sum:prop}
Let $E$ be an elliptic curve over a finite field $\Fq$. Then
\begin{equation}
\label{form_sum:eqn}
\sum_{[E']_\Fq\in\BE(\Fqb/\Fq;E)}{1\over\#\Aut_\Fq(E')} = 1 .
\end{equation}
\end{prop}

\begin{proof}
In the discussion above, take $K=\Fq$ and $L=\Fqb$.  Since 
$\Gal(\Fqb/\Fq)\cong\hat{\BZ}$, the exact sequence~(\ref{seq:eqn}) becomes
\[ 1 \lra A \lra B \mapright{\pi} \hat{\BZ} \lra 1 . \]
Since $\hat{\BZ}$ is freely generated as a profinite group by $1$,
a section $s:\hat{\BZ}\ra B$ is determined by $s(1)$, and every element
of $\pi^{-1}(1)$ gives rise to a section.  Thus, 
$\#S = \#\pi^{-1}(1) = \#A$, 
and dividing equation~(\ref{gen_form_sum:eqn}) by the finite number
$\#A$ yields~(\ref{form_sum:eqn}).
\end{proof}

\begin{cor}
\label{num_curves:cor}
For every prime power $q$,
\[\#' \Big\{ E : \hbox{ $E$ is an elliptic curve over $\Fq$}\Big\} 
                                              \Big/\cong_\Fq \ = q .\]
\end{cor}

\begin{proof}
Let $T$ be the set of elliptic curves over $\Fq$ up to $\Fqb$-isomorphism
and let $U$ be the set of elliptic curves over $\Fq$ up to $\Fq$-isomorphism.
We know that the $j$-invariant provides a bijection between $T$ and $\Fq$,
so $\#T=q$. Also, $U=\bigcup_{[E]_\Fqb\in T}\BE(E)$, so that
\[
\#'U = \sum_{[E']_\Fqb\in T}\quad\sum_{[E]_\Fq\in\BE(E')}{1\over\#\Aut_\Fq(E)}
     = \#T  = q 
\]
as claimed.
\end{proof}

There is a result analogous to Proposition~\ref{form_sum:prop} for
the forms of a triple $(E,P,Q)$.

\begin{prop}
\label{pairform_sum:prop}
Let $E$ be an elliptic curve over a finite field $\Fq$, and let
$P,Q\in E(\Fq)$. Then
\[\sum_{[E',P',Q']_\Fq\in\BE(\Fqb/\Fq;E,P,Q)}
                       {1\over\#\Aut_\Fq(E',P',Q')} = 1 ,\]
where $\Aut_\Fq(E',P',Q')$ denotes the subgroup of $\Aut_\Fq(E')$ 
consisting of those automorphisms that fix $P'$ and $Q'$.
\end{prop}

\begin{proof}
This result follows from making the obvious changes in the proof
of Proposition~\ref{form_sum:prop} and the discussion preceding it.
\end{proof}

\begin{notation}
Suppose $L$ is a Galois extension of a field $K$, $E$ is an elliptic curve
over $K$, and $F$ is an $L/K$-form of $E$.  Given an isomorphism $f:E_L\ra F_L$
and an element $\sigma$ of $\Gal(L/K)$, let $f^\sigma$ be the isomorphism
$(1\times\sigmab)\circ f\circ(1\times\sigmab)^{-1}:E_L\ra F_L$
(here one of the $1\times\sigmab$'s is a $K$-scheme automorphism of $E_L$, 
and the other is a $K$-scheme automorphism of $F_L$).  
If $f$ is defined locally by 
polynomials with coefficients in $L$, then $f^\sigma$ is defined by
the same polynomials with $\sigma$ applied to the coefficients.
\end{notation}

\begin{prop}
\label{automorphism:prop}
Let $E$ be an elliptic curve over a finite field $\Fq$, and let
$\alpha$ be an automorphism of $E_\Fqb$.  Then there is an $\Fqb/\Fq$-form 
$F$ of $E$ and an isomorphism $f:E_\Fqb\ra F_\Fqb$ such that
$\alpha=f^{-1}\circ f^\sigma$, where $\sigma$ is the $q$-th power 
automorphism of $\Fqb$.
\end{prop}

\begin{proof}
With notation as above, let $s:G\ra B$ be the section defined by
sending $\sigma$ to $(\alpha\circ(1\times\sigmab), \sigmab)$
and let $F=E(s)$ and $f=f_s$.
It is not difficult to check that $\alpha=f^{-1}\circ f^\sigma$.
\end{proof}

\section{Modular curves over finite fields}
\label{X:sec}

As indicated in the Introduction, in section~\ref{proof:sec} we will
need to use bounds obtained from modular curves other than the ``standard'' 
modular curves $X(\ell)$ and $X_1(\ell)$.  In this section we define the 
curves we will need, and prove some basic results about them.

First, we recall some facts about Frobenius morphisms of schemes and
elliptic curves (see the discussion in \cite{KaMa}, chapter 12).
For any scheme $S$ over $\Fp$, we define the {\em ($p^r$-th power)
absolute Frobenius} morphism $F_{p^r,\abs}:S\ra S$ to be the morphism 
corresponding to the endomorphism $x\mapsto x^{p^r}$ of affine rings.
If $S$ is a scheme over a field $K$ of characteristic $p>0$, we denote
by $S^{(p^r)}$ the scheme over $K$ defined by the cartesian diagram
% Let's draw a cartesian diagram.
\begin{center}
\begin{picture}(700,700)(-350,-350)
%
% First, the four schemes...
\put(-450, 250){\makebox(0,0){$S^{(p^r)}$}}
\put( 450, 250){\makebox(0,0){$S$}}
\put(-450,-250){\makebox(0,0){$\Spec(K)$}}
\put( 450,-250){\makebox(0,0){$\Spec(K)$}}
%
% Next, the two vertical arrows...
\put(-450, 180){\vector( 0,-1){360}}
\put( 450, 180){\vector( 0,-1){360}}
%
% Then the two horizontal arrows...
\put(-340, 250){\vector( 1, 0){740}}
\put(-250,-250){\vector( 1, 0){500}}
%
% And finally, the names of the horizontal morphism.
\put(   0,-180){\makebox(0,0){$F_{p^r,\abs}$}}
\end{picture}
\end{center}
so that if $S$ is defined locally by polynomials $f_i\in K[x_1,\ldots,x_n]$
then $S^{(p^r)}$ is defined locally by the polynomials $f_i^{(p^r)}$ 
obtained from the $f_i$ by raising all the coefficients to the $p^r$-th
power.

In view of the cartesian property of the above diagram, the $p^r$-th power
absolute Frobenius on $S$ factors through $S^{(p^r)}$; that is, there is a 
morphism 
$F_{p^r}=F_{p^r,S/K}:S\ra S^{(p^r)}$ of $K$-schemes, called
the {\em ($p^r$-th power relative-to-$K$) Frobenius}, such that 
$F_{p^r}$ composed with the map from $S^{(p^r)}$ to $S$ is the 
morphism $F_{p^r,\abs}$ on $S$.  If $S$ is affine and defined by polynomials
$f_i$ as above, then $F_{p^r}$ takes a point $P=(a_1,\ldots,a_n)$ on $S$
to the point $P^{(p^r)}=(a_1^{p^r},\ldots,a_n^{p^r})$ on $S^{(p^r)}$.
In the special case where $S$ is an elliptic curve $E$ over $K$, there  is 
a natural elliptic curve structure on $E^{(p^r)}$, and the Frobenius
$F_{p^r}$ is actually an isogeny.  The dual isogeny of $F_{p^r}$ (see
\cite{Si}, section III.6) is the {\em Verschiebung} 
$V_{p^r}:E^{(p^r)}\ra E$, and the composed map 
$V_{p^r}\circ F_{p^r} : E\ra E$ is the multiplication-by-$p^r$ map on $E$.

We also recall that an elliptic curve $E$ over a field $K$ of characteristic
$p>0$ is called {\em supersingular} if $E$ has no $\Kb$-defined points
of order $p$ (see \cite{Si}, section V.3).  This is equivalent
to the condition that for some $r>0$ the only $\Kb$-valued point in the
kernel of the Verschiebung $V_{p^r}$ is the zero point (which implies the
same statement for {\em all} $r>0$).

The following notation will be useful in this section and the next.

\begin{notation}
Suppose $p$ is a prime number and $m$ and $n$ are positive integers with 
$m | n$ and $m$ coprime to $p$, and write $n=n'p^r$ with $n'$ coprime to $p$.
If $K$ is a field of characteristic $p$ containing a primitive $m$-th root 
of unity $\zeta_m$ and $L$ is an extension field of $K$, we denote by
$Z(L/K; \zeta_m, m, n)$ the set of $L$-isomorphism classes
\[ Z(L/K; \zeta_m, m, n) = \Big\{ (E,P,Q,R) : 
   \vtop{\hsize=4in\raggedright\noindent
   $E$ is an elliptic curve over $K$,
   $P,Q \in E(K)$ with $\ord P = m$ and $\ord Q = n'$
   and $e_m(P, (n'/m) Q) = \zeta_m$, and
   $R \in E^{(p^r)}(K)$ such that $R_\Kb$ generates the kernel of
   the Verschiebung $V_{p^r}: E^{(p^r)}_\Kb \ra E_\Kb \Big\} \Big/\cong_L$ }\]
where $\ord P$ is the order of $P$ in the group $E(K)$ and $e_m$ is the 
Weil pairing on $E[m]$ (see \cite{Si}, section III.8), and
where two such quadruples $(E,P,Q,R)$ and
$(E',P',Q',R')$ are said to be $L$-isomorphic if there is a
$L$-isomorphism $f:E_L\ra E'_L$ such that $f$ takes $P_L$ to $P'_L$
and $Q_L$ to $Q'_L$ and such that $f^{(p^r)}$ takes $R_L$ 
to $R'_L$. Denote by $[E,P,Q,R]_L$ the $L$-isomorphism class of the
quadruple $(E,P,Q,R)$.

Also, we denote by $Y(L/K;\zeta_m, m, n)$ the set of $L$-isomorphism classes
\[ Y(L/K; \zeta_m, m, n) = \Big\{ (E,P,Q) : 
   \vtop{\hsize=4in\raggedright\noindent
   $E$ is an elliptic curve over $K$,
   $P,Q \in E(K)$ with $\ord P = m$ and $\ord Q = n$
   and $e_m(P, (n/m) Q) = \zeta_m \Big\} \Big/\cong_L$ }\]
where two such triples $(E,P,Q)$ and $(E',P',Q')$ are said to be
$L$-isomorphic if there is an $L$-isomorphism $f:E_L\ra E'_L$
that takes $P_L$ to $P'_L$ and $Q_L$ to $Q'_L$. Denote by $[E,P,Q]_L$
the $L$-isomorphism class of the triple $(E,P,Q)$.
\end{notation}

\begin{prop}
\label{existence:prop}
Let $q=p^e$ be a prime power, suppose $m$ and $n$ are positive integers
such that $m | \gcd(n, q-1)$, write $n=n'p^r$ with $n'$ coprime to $p$,
and pick a primitive $m$-th root of unity $\zeta_m \in\Fq$.
There exists a proper nonsingular irreducible curve
$\Xb(m,n)$ over $\Fqb$ provided with a map 
$J:\Xb(m,n)\ra\BP^1_\Fqb\supset\BA^1_\Fqb=\Spec(\Fqb[j])$
with the following properties:
\begin{enumerate}
\item 
\label{cor:item}
There is a natural bijection between the set of finite points of $\Xb(m,n)$ 
(that is, the points in $J^{-1}(\BA^1)$) and the set 
$Z(\Fqb/\Fqb; \zeta_m, m,n)$.

\item
\label{jcor:item}
The bijection given in \ref{cor:item} has the property
that if $x\in\Xb(m,n)$ corresponds to 
$[E,P,Q,R]_\Fqb$ then $J(x)=j(E)$, the $j$-invariant of $E$.

\item
\label{Fqdef:item}
$\Xb(m,n)$ can be defined naturally over $\Fq$; that is, there is a 
proper nonsingular irreducible curve $X_q(m,n)$ over $\Fq$
and an isomorphism
\begin{equation}
\label{Fqdef:eqn}
\Xb(m,n)\cong X_q(m,n)\times_{\Spec(\Fq)}\Spec(\Fqb)
\end{equation}
such that
the  $q$-th power relative-to-$\Fqb$ Frobenius map $F:\Xb(m,n)\ra \Xb(m,n)$ 
obtained from the isomorphism~(\ref{Fqdef:eqn}) and the canonical 
identification 
\[\left(X_q(m,n)\times_{\Spec(\Fq)}\Spec(\Fqb)\right) = 
         \left(X_q(m,n)\times_{\Spec(\Fq)}\Spec(\Fqb)\right)^{(q)}\]
has the property that if the point $x\in \Xb(m,n)$ corresponds to 
$[E,P,Q,R]_\Fqb$, then $F(x)$ corresponds to
$[E^{(q)},P^{(q)},Q^{(q)},R^{(q)}]_\Fqb$.

\end{enumerate}
\end{prop}

\begin{proof}
We will rely heavily on results from~\cite{KaMa}.

First consider the case where $n'>2$.

Pick a primitive $n'$-th root of unity $\zeta_{n'}\in \Fqb$ such that
$\zeta_m=\zeta_{n'}^{n'/m}$, let $\Xb(n',n)$ be the 
$\Fqb$-scheme denoted in \cite{KaMa} by
$\Mb([\Gamma(n')]^\can,[\Ig(p^r)])$ (in~\cite{KaMa}, see sections~4.3 
and~8.6 for the definition of $\Mb(\cdot)$, sections~3.1 and~9.1
for the definition of $[\Gamma(n')]^\can$, and section~12.3 for the 
definition of $[\Ig(p^r)]$), and let
$J':\Xb(n',n)\ra\BP^1_\Fqb\supset\BA^1_\Fqb=\Spec(\Fqb[j])$ be the natural
map from $\Xb(n',n)$ to the ``$j$-line''~$\BP^1_\Fqb$ defined in
section~8.2 of~\cite{KaMa}.
By their very definitions, $\Xb(n',n)$ and $J'$ satisfy
statements~\ref{cor:item} and~\ref{jcor:item} of the proposition
(with $m$ replaced by $n'$ and $J$ replaced by $J'$), and from
Corollary~12.7.2 (page~368) of~\cite{KaMa},
whose hypotheses are satisfied when $n'>2$,
we see that $\Xb(n',n)$ is
a proper nonsingular irreducible curve.
{}From chapter~7 of~\cite{KaMa}, we know the group
\[G=\Big(\SL_2(\BZ/n'\BZ) \times (\BZ/p^r\BZ)^*\Big)\Big/\pm 1\]
(where the group $\{\pm 1\}$ is embedded diagonally in the product)
acts on the covering $\Xb(n',n)$ of $\BP^1$;
the action is such that an element
\[ \pm \left(\left(\matrix{a&b\cr c&d\cr}\right), u \right) \]
of $G$ takes the point corresponding to the class 
$[E,P,Q,R]_\Fqb\in Z(\Fqb/\Fqb; n',n)$
to the point corresponding to the class $[E,aP+cQ, bP+dQ, uR]_\Fqb$.
In fact, from Corollaries~10.13.12 (page~336) and~12.9.4 
(page~381) of~\cite{KaMa} we see that the degree of $\Xb(n',n)$ over $\BP^1$
is equal to $\#G$; since $G$ acts faithfully on $\Xb(n',n)$, this shows
that $\Xb(n',n)$ is a Galois covering of $\BP^1$ with group $G$.

Define a subgroup $H$ of $G$ by
\[ H = \left\{ \pm \left(\left(\matrix{1&0\cr a&1\cr}\right), 1\right)\in G :
       a \equiv 0 \bmod{m}\right\} \]
and define $\Xb(m,n)$ to be the quotient of $\Xb(n',n)$ by the group $H$.
Let 
$J:\Xb(m,n)\ra\BP^1$ be the map induced from $J'$.

Now, a finite point on $\Xb(m,n)$ corresponds to an $H$-orbit of
the finite points on $\Xb(n',n)$; thus, the finite points on $\Xb(m,n)$
correspond to the $\Fqb$-isomorphism classes of sets of the form
\[\Big\{ (E, P+aQ, Q, R) : a\equiv 0 \bmod{m}\Big\},\]
where $[E,P,Q,R]_\Fqb\in Z(\Fqb/\Fqb; n',n)$ and where
two such sets $\{(E,P+aQ,Q,R)\}$ and $\{(E',P'+aQ',Q',R')\}$
are $\Fqb$-isomorphic if there is an isomorphism $f:E\ra E'$ such that
$f$ maps $Q$ to $Q'$ and the set $\{P+aQ\}$ to $\{P'+aQ'\}$ and 
such that $f^{(p^r)}$ maps $R$ to $R'$.
But there is a natural bijection between the set of all such 
$\Fqb$-isomorphism classes and the set $Z(\Fqb/\Fqb; m,n)$
given by sending the class of $\{(E,P+aQ,Q,R)\}$ to the class
$[E,{n'\over m}P,Q,R]_\Fqb$.  Thus,
$\Xb(m,n)$ satisfies the property given in statement~\ref{cor:item}
of the proposition.

That $J$ satisfies the property given in statement~\ref{jcor:item}
is a consequence of the fact that $J'$ satisfies the corresponding property
and of the construction just given.

Finally, that $\Xb(m,n)$ may be defined over $\Fq$ in the manner
described in statement~\ref{Fqdef:item} follows from general principles
given in~\cite{KaMa} (see in particular the discussion in section~12.10)
and from the fact that the
correspondence in statement~\ref{cor:item} refers only to structures
(in particular, the element $\zeta_m$) that are defined over $\Fq$.

This completes the proof for the case where $n'>2$.  Now suppose
$n'\le 2$.  
The problem with proceeding exactly as before is that the 
results in~\cite{KaMa} that we used in the case $n'>2$ (in particular,
Corollaries~12.7.2,~10.13.12 and~12.9.4) don't apply when $n'\le 2$, because,
in the language of \cite{KaMa}, $[\Gamma(n')]^\can$ is not representable
when $n'\le 2$.
Thus, we have to make some very minor modifications to our previous argument,
although the general idea is exactly the same.  

If $n'=2$ let $f=2$; if $n'=1$ and $p\ne 3$ let $f=3$; if $n'=1$ and
$p=3$ let $f=4$. Consider the
curve $\Xb(fn', fn)$, which, as before, is a Galois covering
of $\BP^1$ with Galois group
\[G=\Big(\SL_2(\BZ/fn'\BZ) \times (\BZ/p^r\BZ)^*\Big)\Big/\pm 1,\]
and which has an interpretation as in statement~\ref{cor:item}.
Now let $H$ be the subgroup
\[H=\left\{\pm\left(\left(\matrix{a&b\cr c&d\cr}\right), 1 \right)\in G :
          a\equiv d\equiv 1 \bmod{{n'}}, b\equiv 0 \bmod{{n'}},
          \mbox{\ and\ } c\equiv 0 \bmod {m} \right\}, \]
and let $\Xb(m,n)$ be the quotient of $\Xb(fn',fn)$ by $H$.
The proof follows exactly as before.

Thus, the proposition is valid for all values of $n'$.
\end{proof}

There are two special kinds of points on the curves $\Xb(m,n)$
that we will need to keep track of.

\begin{defn}
Let $q$, $m$, $n$, $\Xb(m,n)$, and $J$ be as in 
Proposition~\ref{existence:prop}.  A point $x\in\Xb(m,n)$ is a
{\em cusp} if $x$ is an element of $J^{-1}(\infty)$.  A point of 
$\Xb(m,n)$ which is not a cusp is called a {\em finite point}.  A finite 
point of $\Xb(m,n)$ is a {\em supersingular point} if it corresponds 
to an equivalence class $[E,P,Q,R]_\Fqb$ with a supersingular $E$.
\end{defn}

\begin{notation}
We denote by $g_q(m,n)$ the genus of $\Xb(m,n)$, by $c_q(m,n)$ the number
of cusps of $\Xb(m,n)$, and by $s_q(m,n)$ the number of supersingular
points of $\Xb(m,n)$.
\end{notation}

\begin{prop}
\label{bounds:prop}
For all $q=p^e$, $m$, and $n=n'p^r$ as in Proposition~\ref{existence:prop}
we have
\begin{equation}
\label{g_bound:eqn}
g_q(m,n)\le {1\over 24} m \varphi(n)\psi(n)
\end{equation}
\begin{equation}
\label{c_bound:eqn}
c_q(m,n)\le  \varphi(n)\psi(n)
\end{equation}
and when $p|n$ (that is, when $r>0$) we have
\begin{equation}
\label{s_bound:eqn}
s_q(m,n)\le {1\over 3} m \varphi(n)\psi(n).
\end{equation}
\end{prop}

\begin{proof}
As in the preceding proof, we first assume that $n'>2$.

Let the groups $G$ and $H$ be as in the proof of 
Proposition~\ref{existence:prop}, so that $\Xb(m,n)$ is the quotient
of $\Xb(n',n)$ by $H$.  
{}From Corollary~10.13.12 (page~336) and Corollary~12.9.4 (page~381)
of \cite{KaMa} we find that
\begin{equation}
\label{exact_genus:eqn}
g_q(n',n) =\cases{\displaystyle 1 + {1\over 24}(n-6)\varphi(n)\psi(n)
						& if $n'=n$;\cr
		    &\cr
		    \displaystyle 1 + {1\over 48}(n-12)\varphi(n)\psi(n')
						& if $n'<n$\cr} 
\end{equation}
\[c_q(n',n) = {1\over 2}\varphi(n)\psi(n')\]
\begin{equation}
\label{exact_ss:eqn}
s_q(n',n) = {{p-1}\over 24} n' \varphi(n')\psi(n').
\end{equation}
Since $\#H = n'/m$, the Riemann-Hurwitz formula (\cite{Si}, Theorem~5.9,
page~41) gives us the estimate
\[g_q(m,n) \le\cases{\displaystyle 1+{1\over 24}{m\over n}(n-6)
						\varphi(n)\psi(n)
						& if $n'=n$;\cr
	     &\cr
	          \displaystyle 1+{1\over 48}{m\over n'}(n-12)
						\varphi(n)\psi(n')
						& if $n'<n$,\cr} \]
which leads to~(\ref{g_bound:eqn}).

We also have the trivial bound
\[c_q(m,n)\le c_q(n',n)=
         {1\over 2}\varphi(n)\psi(n')\le{1\over 2}\varphi(n)\psi(n),\]
which certainly implies~(\ref{c_bound:eqn}).

To get a good bound for $s_q(m,n)$, we need to determine necessary
conditions for an element  of $H$ to fix a finite point of $\Xb(n',n)$.
So suppose $x$ is a finite point
of $\Xb(n',n)$, corresponding to the class $[E,P,Q,R]_\Fqb$;
for a non-trivial element of $H$ to fix $x$, we must have
$[E,P,Q,R]_\Fqb=[E,P+aQ,Q,R]_\Fqb$ for some 
$a$ with $a\equiv 0\bmod{m}$ and $a\not\equiv 0\bmod{n'}$,
so there must be an automorphism $\alpha$ of $E$
that fixes $Q$ and sends $P$ to $P+aQ$. Thus $\alpha\ne\pm 1$, and from
Corollary~2.7.1 (page~85) of~\cite{KaMa} we see that $\alpha$ satisfies
$\alpha^2-t\alpha+1=0$ for some integer $t$ with $|t|\le 1$.  In
particular, this means that $(2-t)Q = 0$, which is impossible if
$n'>3$.  Thus, if $n'>3$ no non-trivial element of $H$ fixes any finite
point of $\Xb(n',n)$, so every finite point of $\Xb(m,n)$ has
$\#H$ points of $\Xb(n',n)$ lying over it; this gives us
\[s_q(m,n)={m\over n'} s_q(n',n) = {{p-1}\over 24}m\varphi(n')\psi(n').\]
When $n'=3$, we at least have the bound
\[s_q(m,n)\le s_q(n',n)\le {{p-1}\over 8}m\varphi(n')\psi(n'),\]
so that in any case if $p| n$ we have
\[s_q(m,n)\le{1\over 24}m\varphi(n)\psi(n).\]
This gives us~(\ref{s_bound:eqn}).

Thus, when $n'>2$, the inequalities of the proposition hold.

When $n'\le 2$, let $f$, $G$, and $H$ be
as in the case $n'\le 2$ of the proof of Proposition~\ref{existence:prop},
so that $\Xb(m,n)$ is the quotient of $\Xb(fn',fn)$ by $H$.
Once again, one can check that equation~(\ref{exact_genus:eqn}) 
and the Riemann-Hurwitz
formula lead to~(\ref{g_bound:eqn}).

To prove~(\ref{c_bound:eqn}), we note that 
it is possible to define $\Xb(m,n')$ as the quotient of $\Xb(fn',fn)$ by
the subgroup of $G$ generated by $H$ and the image of $(\BZ/p^r\BZ)^*$
in $G$; this gives us a map from $\Xb(m,n)$ to $\Xb(m,n')$ consistent
with the maps from these curves to $\BP^1$ and of degree at most
$\varphi(p^r)$, so that $c_q(m,n)\le \varphi(p^r) c_q(m,n')$.  From this
inequality we see that it suffices to prove~(\ref{c_bound:eqn}) when
$n=n'$, that is, when $r=0$.  But from statement~1 of Theorem~10.9.1 
(page~301) of~\cite{KaMa} we can calculate that
$c_q(2,2)=3$, $c_q(1,2)=2$, and $c_q(1,1)=1$, so 
inequality~(\ref{c_bound:eqn}) {\em does} hold when $r=0$.

Finally, suppose $p|n$.  Using the trivial bound $s_q(m,n)\le s_q(fn',fn)$
and equation~(\ref{exact_ss:eqn}), we see that
\[ s_q(m,n)\le {p-1 \over 24} fn'\varphi(fn')\psi(fn'); \]
it is easy to check that this inequality implies~(\ref{s_bound:eqn}),
except when $n=p=3$.  But in this case we notice
that $G=H$, so that $\Xb(1,3)=\BP^1$ has exactly one supersingular
point (corresponding to the elliptic curve with $j$-invariant 0),
and we can verify~(\ref{s_bound:eqn}) directly.

Thus, inequalities~(\ref{g_bound:eqn}),~(\ref{c_bound:eqn}), 
and~(\ref{s_bound:eqn}) hold in every case.
\end{proof}

\begin{remark}
{}From equation~(\ref{exact_genus:eqn}) we see that $1/24$ is the 
smallest possible constant in inequality~(\ref{g_bound:eqn}).
The facts that $c_q(1,1)=1$ and  $s_2(1,2)=1$ show that equality
is sometimes obtained in inequalities~(\ref{c_bound:eqn}) 
and~(\ref{s_bound:eqn}).
\end{remark}

We now focus on the curves $X_q(m,n)$.  In particular, we may ask whether
there is a modular interpretation for the $\Fq$-defined points of 
$X_q(m,n)$.  The answer is ``yes''.

\begin{prop}
\label{Fqdef:prop}
Let $q$, $m$, $n=n'p^r$, $\zeta_m$, and $X_q(m,n)$ be as in 
Proposition~\ref{existence:prop}.  
There is a bijection between the set of finite points of $X_q(m,n)(\Fq)$ 
(that is, the finite points of $X_q(m,n)$ that are defined over $\Fq$) 
and the set $Z(\Fqb/\Fq; \zeta_m, m, n)$.
\end{prop}

\begin{proof}
Let $F:\Xb(m,n)\ra \Xb(m,n)$ be the $q$-th power relative-to-$\Fqb$ Frobenius 
map, as in statement~\ref{Fqdef:item} of Proposition~\ref{existence:prop}. 
Then there is a bijection between $X_q(m,n)(\Fq)$ and the set of points
of $\Xb(m,n)$ fixed by $F$, given by $x\mapsto x_\Fqb$.
Again by statement~\ref{Fqdef:item} of Proposition~\ref{existence:prop}, 
we know that the finite points of this last set correspond to
the elements of the set
\[ S=\Big\{ [E,P,Q,R]_\Fqb\in Z(\Fqb/\Fqb; \zeta_m, m, n) :
        [E,P,Q,R]_\Fqb = [E^{(q)},P^{(q)},Q^{(q)},R^{(q)}]_\Fqb \Big\} . \]
Thus, we need only show that there is a bijection between
the sets $S$ and $Z(\Fqb/\Fq; \zeta_m, m, n)$.

There is clearly an injective map from $Z(\Fqb/\Fq; \zeta_m, m, n)$ to 
$S$ defined by 
sending  $[E,P,Q,R]_\Fqb$ to $[E_\Fqb,P_\Fqb,Q_\Fqb,R_\Fqb]_\Fqb$.
We need only show that this map is surjective.

Suppose $[E,P,Q,R]_\Fqb$ is an element of $S$, and let
$f:E\ra E^{(q)}$ be an isomorphism that takes the quadruple
$(E,P,Q,R)$ to the quadruple $(E^{(q)},P^{(q)},Q^{(q)},R^{(q)})$.
Since $E\cong_\Fqb E^{(q)}$,
we have $j(E)=j(E^{(q)})=(j(E))^q$, so $j(E)\in\Fq$.  Let $E'$ be any 
elliptic curve over $\Fq$ with $j(E')=j(E)$; since elliptic curves over
$\Fqb$ are classified up to $\Fqb$-isomorphism by their $j$-invariants, there 
is an isomorphism $g:E\ra E'_\Fqb$.  By Proposition~\ref{automorphism:prop},
there is a form $F$ of $E'$ and an isomorphism $h:E'_\Fqb\ra F_\Fqb$
such that 
\[ g\circ f^{-1}\circ (g^{(q)})^{-1} = h^{-1}\circ h^{(q)} ,\]
and by replacing $E'$ with $F$ and $g$ with $h\circ g$,
we may assume that $g^{(q)}\circ f\circ g^{-1}$ is the identity
on $E'_\Fqb$.

Now consider the point $g(P)\in E'_\Fqb(\Fqb)$.  We  have
\[ g(P) = (g^{(q)}\circ f\circ g^{-1})(g(P)) = g^{(q)}(f(P)) =
          g^{(q)}(P^{(q)}) = (g(P))^{(q)} , \]
so $g(P)$ is an $\Fq$-defined point of $E'_\Fqb$; that is, there is a 
point $P'\in E'(\Fq)$ such that $g(P)=P'_\Fqb$.  Similarly, we see
that $g(Q)$ and $g^{(p^r)}(R)$ come from points $Q'\in E'(\Fq)$ and
$R'\in E'^{(p^r)}(\Fq)$, so that $[E',P',Q',R']_\Fqb$ is an element
of $Z(\Fqb/\Fq; \zeta_m, m, n)$ that maps to the element 
$[E,P,Q,R]_\Fqb$ of $S$.  Thus,
the natural map from $Z(\Fqb/\Fq; \zeta_m, m, n)$ to $S$ is bijective, 
and the proposition is proven.
\end{proof}

\begin{remark}
More generally, if $K$ is any field containing 
$\Fq$ and $\Kb$ is the algebraic closure
of $K$, we know from Lemma~8.1.3.1 (page~225) of~\cite{KaMa} that 
there is a bijection between the set of finite 
$\Kb$-valued points of $X_q(m,n)$ and $Z(\Kb/\Kb;\zeta_m,m,n)$, and
we may ask whether the finite $K$-valued points of
$X_q(m,n)$ correspond to the elements of $Z(\Kb/K;\zeta_m,m,n)$.
The proof of Proposition~3.2 (page~274) of~\cite{DeRa} provides
an answer:  The obstruction to a $K$-valued point giving rise to
a quadruple $(E,P,Q,R)$ defined over $K$ lies in a certain $H^2$, and
it is shown in the proof of Proposition~3.2 of~\cite{DeRa} that
this obstruction is zero.  In the special case $K=\Fq$ we consider
above, the argument simplifies, because in this case the whole
$H^2$ where the obstruction lives is trivial.  One can use this
argument to provide a more conceptual proof of 
Proposition~\ref{Fqdef:prop}.  The interested reader should 
consult~\cite{DeRa}.
\end{remark}

\begin{cor}
\label{y_estimate:cor}
There is a constant $C'\le1/12+5\sqrt{2}/6\approx 1.262$
such that for all $q$, $m$, and $n=n'p^r$ as in 
Proposition~\ref{existence:prop}
the following statements are true:
\begin{enumerate}
\item
\label{non_Ig:item}
If $n'=n$, then there is a bijection between the set 
$Y(\Fqb/\Fq; \zeta_m, m, n)$ and the set of finite points of $X_q(m,n)(\Fq)$.

\item
\label{Ig:item}
If $n'<n$, then there is a bijection between the set 
$Y(\Fqb/\Fq; \zeta_m, m, n)$ and the 
set of finite non-supersingular points of $X_q(m,n)(\Fq)$.

\item
\label{y_estimate:item}
We have the estimate
\begin{equation}
\label{Yerror:eqn}
\left|\#Y(\Fqb/\Fq;\zeta_m,m,n) - q \right| 
                 \le C' m \varphi(n)\psi(n)\sqrt{q}.
\end{equation}
\end{enumerate}
\end{cor}

\begin{proof}
If $n'=n$ then there is a bijection between the sets 
$Y(\Fqb/\Fq;\zeta_m, m, n)$ and
$Z(\Fqb/\Fq; \zeta_m, m, n)$, given by mapping 
$[E,P,Q]_\Fqb$ to 
$[E,P,Q,O]_\Fqb$, where $O$ is the zero element of $E=E^{(1)}$ (which generates
the kernel of the Verschiebung $V_1$, the identity map).
Thus, statement~\ref{non_Ig:item} follows immediately from
Proposition~\ref{Fqdef:prop}.

If $n'<n$, let 
\[Z'(\Fqb/\Fq;\zeta_m,m,n)=\left\{[E,P,Q,R]_\Fqb\in Z(\Fqb/\Fq;\zeta_m,m,n) : 
   \mbox{$E$ is not supersingular}\right\}.\]
Let $M$ be the map from $Y(\Fqb/\Fq;\zeta_m,m,n)$ to 
$Z(\Fqb/\Fq;\zeta_m,m,n)$ that sends 
$[E,P,Q]_\Fqb$ to  $[E,P,p^rQ,(n'Q)^{(p^r)}]_\Fqb$. 
The image of $M$ lies in $Z'(\Fqb/\Fq;\zeta_m,m,n)$, 
because if $Q\in E(\Fq)$ has order $n$ then $n'Q$ has order $p^r\ne1$,
so that $E$ is not supersingular.  Choose integers $a$ and $b$ such that
$ap^r+bn'=1$; 
then the inverse of $M$ is 
the map from $Z'(\Fqb/\Fq;\zeta_m,m,n)$ to 
$Y(\Fqb/\Fq;\zeta_m,m,n)$ 
that sends $[E,P,Q,R]_\Fqb$ to $[E,P,(aQ+bR')]_\Fqb$, where $R'$
is the element of $E(\Fq)$ such that $(R')^{(p^r)}=R$ --- this element 
exists and is unique because $\Fq$ is perfect.
Thus $M$ is a bijection between $Y(\Fqb/\Fq;\zeta_m,m,n)$ and 
$Z'(\Fqb/\Fq;\zeta_m,m,n)$, so that 
statement~\ref{Ig:item} follows from Proposition~\ref{Fqdef:prop}.

To prove statement~\ref{y_estimate:item} we will need to use the
Weil conjectures for curves (see~\cite{We:CA} or~\cite{Bo}); 
in particular, we will need the inequality (\cite{We:CA}, Corollaire~3, page~70)
\begin{equation}
\label{Weil_conj:eqn}
\Bigl|\#X_q(m,n)(\Fq)-1-q\Bigr| \le 2g_q(m,n)\sqrt{q}.
\end{equation}

First suppose that $n'=n$.  Then statement~\ref{non_Ig:item},
combined with the inequalities~(\ref{g_bound:eqn}),~(\ref{c_bound:eqn}),
and~(\ref{Weil_conj:eqn}), gives us
\[ \left|\#Y(\Fqb/\Fq;\zeta_m,m,n) - q\right| \le 
            1 + \varphi(n)\psi(n) + {1\over 12}m\varphi(n)\psi(n)\sqrt{q}.\]
On the other hand, if $n'<n$, then statement~\ref{Ig:item}, combined
with the 
inequalities~(\ref{g_bound:eqn}),~(\ref{c_bound:eqn}),~(\ref{s_bound:eqn}),
and~(\ref{Weil_conj:eqn}), gives us
\[ \left|\#Y(\Fqb/\Fq;\zeta_m,m,n) - q\right| \le 
            1 + \varphi(n)\psi(n) + {1\over 3}m\varphi(n)\psi(n)
                    + {1\over 12}m\varphi(n)\psi(n)\sqrt{q}. \]
Thus, statement~\ref{y_estimate:item} will hold if we choose $C'$ so
that for all $q$, $m$, and $n$ we have
\[ C' \ge {1\over m\varphi(n)\psi(n)\sqrt{q}}+{1\over m\sqrt{q}}
          + {1\over 3\sqrt{q}} + {1\over 12} . \]
However, since $\#Y(\Fqb/\Fq;1,1,1)=q$ (as we noted in the proof 
of Corollary~\ref{num_curves:cor}, where the set was called $T$), 
we need only have the above
inequality when $n>1$.  Thus, $C'=1/12+5\sqrt{2}/6$ will do.
\end{proof}

With inequality~(\ref{Yerror:eqn}) in hand, we can proceed to the calculations
of section~\ref{proof:sec}.

\section{Proof of the theorem}
\label{proof:sec}

Fix a prime power $q=p^e$, and let $\zeta_{q-1}$ be a primitive
$(q-1)$-th root of unity in $\Fq$.  For each $m$ dividing $q-1$, let
$\zeta_m$ be the primitive $m$-th root of unity $\zeta_{q-1}^{(q-1)/m}$.
Recall that for every pair $(m,n)$ of positive integers with
$m$ dividing $\nb$ we have sets
$Y(\Fq/\Fq; \zeta_m, m,n)$ and $Y(\Fqb/\Fq; \zeta_m, m,n)$.
For each pair $(m,n)$ with $m|n$ we also define a set
\[ W(\Fq;m,n)=\Big\{ E/\Fq : E[n](\Fq)\cong (\BZ/m\BZ)\times(\BZ/n\BZ)\Big\}
                                                           \Big/\cong_\Fq .\]
Note that $W(\Fq;m,n)$ is empty unless $m$ divides $q-1$; 
see Corollary~8.1.1 (page~98) of~\cite{Si}.
Also, for every positive integer $N$, we have the set $V(\Fq; N)$.
Our goal is to estimate the weighted cardinality of $V(\Fq; N)$.

For all the appropriate values of $m$, $n$, and $N$, let
$v(N)=\#'V(\Fq;N)$ and $w(m,n)=\#'W(\Fq;m,n)$ and 
$y(m,n)=\#Y(\Fqb/\Fq; \zeta_m, m, n)$
(note that $y(m,n)$ is a non-weighted cardinality).
Corollary~\ref{y_estimate:cor} gives us an estimate 
for $y(m,n)$ for all pairs $(m,n)$ with $m$ dividing $\nb$.  To get
{}from these estimates to an estimate for $v(N)$, we need to make 
explicit the relationships among the sets mentioned above.

\begin{notation}
Let $t$ and $u$ denote the multiplicative arithmetic functions defined 
on prime powers $\ell^a$ by $t(\ell^a)=\ell^{\lfloor a/2 \rfloor}$ and
$u(\ell^a)=\ell^{\lceil a/2 \rceil}$; thus, for every positive
integer $N$ we have $N/t(N)^2 = u(N)^2/N$, and this number is a
squarefree integer.
Also, given a positive integer $n$ and a prime number $\ell$, we will 
denote by $n_{(\ell)}$ the largest power of $\ell$ dividing $n$. 
Thus, for example, $t(24)=2$ and $u(24)=12$ and $24_{(2)} = 8$.
\end{notation}

\begin{lemma}
\label{v_from_w:lemma}
Let $N$ be any positive integer.  Then
\begin{equation}
\label{VfromW:eqn}
V(\Fq;N) = \coprod_{d | {\gcd(u(N),q-1)}} W \left( \Fq; d , 
                                    {N\over \gcd(d, t(N))} \right), 
\end{equation}
and
\begin{equation}
\label{v_from_w:eqn}
v(N) = \sum_{d | {\gcd(u(N),q-1)}} w \left( d , 
                                    {N\over \gcd(d, t(N))} \right).
\end{equation}
\end{lemma}

\begin{proof}
Since (\ref{v_from_w:eqn}) follows from (\ref{VfromW:eqn}), it suffices to
prove (\ref{VfromW:eqn}).  Also, (\ref{VfromW:eqn}) is equivalent to
\begin{equation}
\label{VfromWunr:eqn}
V(\Fq;N) = \coprod_{d | u(N)} W \left( \Fq; d , 
                                    {N\over \gcd(d, t(N))} \right), 
\end{equation}
because the additional sets we get in (\ref{VfromWunr:eqn}) are all empty.

It is easy to see that $W(\Fq;d,N/\gcd(d,t(N)))\subset V(\Fq;N)$ 
for each divisor $d$ of $u(N)$.
On the other hand,
suppose we are given an elliptic curve $E$ over $\Fq$ with 
$[E]_\Fq\in V(\Fq;N)$.  It is not hard to show that
if $d|u(N)$ then $[E]_\Fq$ is an element of  $W(\Fq;d,N/\gcd(d,t(N)))$
if and only if $d$ is the largest divisor of $u(N)$ for which
$\#E[d](\Fq)=d^2$; this is easy to check when $N$ is a prime power,
and it suffices to check only this case because
for all pairs $(m,n)$ with $m|n$ we have 
\[W(\Fq;m,n)=\bigcap_{\mbox{\footnotesize primes $\ell$}}
                            W(\Fq;m_{(\ell)},n_{(\ell)}).\]
Thus, for every element $[E]_\Fq$ of $V(\Fq;N)$ there is a unique 
divisor $d$ of $u(N)$
with $[E]_\Fq\in W(\Fq;d,N/\gcd(d,t(N)))$, and we are done.
\end{proof}

\begin{lemma}
\label{y_from_w:lemma}
For every pair $(m,n)$ of positive integers with $m$ dividing $\nb$,
we have
\begin{equation}
\label{y_from_w:eqn}
y(m,n) = m\varphi(n)\psi(n) \sum_{d: m| d| \nb} {w(d,n)\over\psi(n/d)}.
\end{equation}
\end{lemma}

\begin{proof}
Consider the map from $Y(\Fq/\Fq; \zeta_m, m, n)$ to 
$\coprod_{d : m|d|\nb} W(\Fq;d,n)$
that takes $[E,P,Q]_\Fq$ to $[E]_\Fq$.  This map is clearly surjective.  

Consider an elliptic curve $E$ over $\Fq$ with $[E]_\Fq\in W(\Fq;d,n)$
for some $d$ with $m|d|\gcd(n,q-1)$.
It is not difficult to check that there are exactly 
$m\varphi(n)\psi(n)/\psi(n/d)$ ways of choosing a pair
$(P,Q)$ of points of $E(\Fq)$ with $\ord P = m$, $\ord Q = n$, and
$e_m(P,(n/m)Q)=\zeta_m$.
Two such pairs $(P,Q)$ and $(P',Q')$ satisfy $(E,P,Q)\cong_\Fq(E,P',Q')$ 
if and only if $(P',Q')$ lies in the $\Aut_\Fq(E)$-orbit of $(P,Q)$,
and the size of this orbit is the index 
$[\Aut_\Fq(E) : \Aut_\Fq(E,P,Q)] = \#\Aut_\Fq(E)/\#\Aut_\Fq(E,P,Q)$.
Summing over the various $\Aut_\Fq(E)$-orbits of such pairs, we obtain
\[\sum_{(P,Q)} {\#\Aut_\Fq(E)\over\#\Aut_\Fq(E,P,Q)} = 
                         {m\varphi(n)\psi(n)\over\psi(n/d)}.\]
Dividing by $\#\Aut_\Fq(E)$ and summing over $\Fq$-isomorphism classes
of $E$ we obtain
\[\sum_{[E,P,Q]_\Fq\in Y(\Fq/\Fq;\zeta_m,m,n)} {1\over\#\Aut_\Fq(E,P,Q)} = 
   \sum_{d : m| d|\nb} {m\varphi(n)\psi(n)\over\psi(n/d)}\#'W(\Fq;d,n) . \]
But the sum on the left hand side is
\[ \sum_{[E',P',Q']_\Fqb\in Y(\Fqb/\Fq;\zeta_m,m,n)}\quad
     \sum_{[E,P,Q]_\Fq\in\BE(E',P',Q')} {1\over\#\Aut_\Fq(E,P,Q)}         \]
and by Proposition~\ref{pairform_sum:prop} this double sum is the
cardinality of $Y(\Fqb/\Fq; \zeta_m, m, n)$.  
This gives us~(\ref{y_from_w:eqn}).
\end{proof}

\begin{lemma}
\label{w_from_y:lemma}
For every pair $(m,n)$ of positive integers with $m$ dividing $\nb$,
we have
\begin{equation}
\label{w_from_y:eqn}
w(m,n) = {\psi(n/m)\over m\varphi(n)\psi(n)}\sum_{j|(\nb/m)} 
					{\mu(j)\over j} y(mj,n) .
\end{equation}
\end{lemma}

\begin{proof}
We calculate:
\begin{eqnarray*}
{w(m,n)\over\psi(n/m)} &=& \sum_{d : m| d| \nb} {w(d,n)\over\psi(n/d)}
                   \sum_{j| {(d/m)}} \mu(j)            \\
             &=& \sum_{j|{\nb/m}} \mu(j)
                   \sum_{d : mj| d|\nb} {w(d,n)\over\psi(n/d)}
    \ =\  \sum_{j|{(\nb/m)}}\mu(j){y(mj,n)\over mj\varphi(n)\psi(n)} 
\end{eqnarray*}
where the last equality follows from~(\ref{y_from_w:eqn}). Multiplying
by $\psi(n/m)$ we get~(\ref{w_from_y:eqn}).
\end{proof}

Now we use the approximation that Corollary~\ref{y_estimate:cor} gives
us for $y(m,n)$ to define approximations for $w(m,n)$ and $v(N)$;
namely, for all pairs $(m,n)$ of positive integers with $m$ dividing
$\nb$, we define
\[\wh(m,n)={q\psi(n/m)\over m\varphi(n)\psi(n)}
                    \sum_{j|{(\nb/m)}}{\mu(j)\over j}
          = {q\psi(n/m)\over m\varphi(n)\psi(n)} 
            \prod_{\ell|{(\nb/m)}}\left(1-{1\over\ell}\right)\]
and for all positive integers $N$ we define
\[\vh(N) = \sum_{d|{{\gcd(u(N),q-1)}}} 
                     \wh\left(d,{N\over\gcd(d,t(N))}\right).\]
We see from Lemma~\ref{w_from_y:lemma} and Corollary~\ref{y_estimate:cor}
that
\begin{eqnarray*}
|w(m,n)-\wh(m,n)| &\le& {\psi(n/m)\over m\varphi(n)\psi(n)} 
                        \sum_{j|{(\nb/m)}}{|\mu(j)|\over j}
                           C' mj\varphi(n)\psi(n)\sqrt{q}  \\
&\ &\\
&\le& C'\psi(n/m)\sqrt{q}\sum_{j|{(\nb/m)}} |\mu(j)|
            \ \le\ C'\psi(n/m)2^{\nu(n)}\sqrt{q}, 
\end{eqnarray*}
where $\nu(n)$ denotes the number of prime divisors of $n$.  From this error
estimate and from Lemma~\ref{v_from_w:lemma}, we find that
\begin{eqnarray}
|v(N)-\vh(N)| &\le& \sum_{d|{{\gcd(u(N),q-1)}}}
                      C'{\psi(N/ d)}2^{\nu(N)}\sqrt{q}
            \ \le\ C'\psi(N)2^{\nu(N)}\sqrt{q}\sum_{d| N}1/d \nonumber \\
        &<&  C'\psi(N)2^{\nu(N)}\sqrt{q}\prod_{\ell| N}{1\over 1-1/\ell}
                \ =\ C'N\rho(N)2^{\nu(N)}\sqrt{q}  .\label{v_error:eqn} 
\end{eqnarray}

To calculate $\wh(m,n)$ and $\vh(N)$, we note that the definition of
$\wh(m,n)$ shows that the ratio $\wh(m,n)/q$ is multiplicative; that is,
\[{\wh(m,n)\over q} = \prod_\ell {\wh(m_{(\ell)},n_{(\ell)})\over q} .\]
This equality, together with the definition of $\vh(N)$, shows that
$\vh(N)/q$ is a multiplicative function of $N$.
A straightforward (if tedious) verification shows that for 
prime powers $\ell^a$ we have
\begin{equation}
\label{prime_power:eqn}
{\vh(\ell^a)\over q}=\left\{ \begin{array}{ll}
            \displaystyle {1\over \ell^{a-1}(\ell-1)}&
                         \mbox{if $q\not\equiv 1 \bmod{\ell^c};$} \\
     \ & \\
            \displaystyle{{\ell^{b+1}+\ell^b-1}\over{\ell^{a+b-1}(\ell^2-1)}}&
                         \mbox{if $q\equiv 1 \bmod{\ell^c},$} 
                            \end{array}
                     \right.
\end{equation}
where $b=\lfloor{a/2}\rfloor$ and $c=\lceil{a/2}\rceil$.

Inequality~(\ref{v_error:eqn}) and equation~(\ref{prime_power:eqn})
show that Theorem~\ref{theorem:thm} will be true if we take $C$ to
be $C'$ and $r(N)$ to be the ratio $\vh(N)/q$.
 \hbox to 2em{}\nobreak\hfill$\hollowsquare$\par

\end{document}